\documentclass[11pt]{article}

\usepackage[latin1]{inputenc}

\usepackage{amsfonts,amsmath,amssymb,amsthm,graphicx,epsfig,float}
\usepackage[T1]{fontenc}

\usepackage[english,francais]{babel}

{
  \newtheorem{theorem}{Theorem}
  \newtheorem*{theorem*}{Theorem}
  \newtheorem{lemma}{Lemma}

  \newtheorem*{corollaire*}{Corollary}
\newtheorem*{proposition*}{Proposition}
\theoremstyle{remark}
  \newtheorem{remark}{Remark}
}

\newcounter{ex}

\newenvironment{rem*}{
  \noindent\textbf{Remarque. }}{}

\newcommand{\Cc}{\mathbb{C}}

\renewcommand {\epsilon}{\varepsilon}

\renewcommand {\leq}{\leqslant}
\renewcommand {\geq}{\geqslant}

\title{{\bf Ahlfors' currents in higher dimension}}
\author{Henry de Thélin}
\date{}

\begin{document}
\maketitle

%$\Cc, \Rr$%

%% Redefinition Titre
\def\figurename{{Fig.}}%
\def\proofname{Preuve}% for AMS-\LaTeX
\def\contentsname{Sommaire}%
%% Fin

\selectlanguage{english}

\begin{abstract}

We consider a nondegenerate
holomorphic map $f: V \mapsto X$ where $(X, \omega)$ is a compact hermitian manifold of
dimension higher or equal to $k$ and $V$ is an open connected complex
manifold of dimension $k$. In this article we give criteria which permit to construct Ahlfors' currents in $X$.

\end{abstract}


AMS: 32A15, 32U40.

\section*{{\bf Introduction}}
\par

Let $f: V \mapsto X$ be a nondegenerate holomorphic map between an
open connected complex manifold $V$ (non-compact) of
dimension $k$ and a compact hermitian manifold $(X,
\omega)$ of
dimension higher or equal to $k$. We consider an exhaustion function
$\tau$ on $V$. It means that (see \cite{Wu}):

(i) $\tau: V \mapsto [0 ,+ \infty[$ is $C^1$.

(ii) $\tau$ is proper
    (i.e. $\tau^{-1}(\mbox{compact})=\mbox{compact}$).

(iii) There exists $r_0 >0$ such that $\tau$ has only isolated
    critical points in $\tau^{-1}( [r_0 , + \infty[ )$.

In this article we will employ the notation $V(r)= \tau^{-1}([0 , r[)$. 

The first important example is $V= \Cc^k$ and $\tau = \| z \|^2$. When $k=1$ we are studying entire
cuves in $X$. An another example is for $V$ a pseudoconvex domain in $\Cc^k$. If
$\tau_0$ is its exhaustion function, we can easily transform $\tau_0$
into a function $\tau$ which satisfies the previous hypothesis (see
\cite{Ra} p. 63-65).

The goal of this article is to construct Ahlfors' currents in $X$ with
$V$ and $f$. By definition, a Ahlfors' current is a limit of a sequence like
$\frac{f_{*}[V(r_n)]}{\mbox{volume}(f(V(r_n)))}$ which is a {\bf closed}
positive current of bidimension $(k,k)$ (here $r_n \rightarrow +
\infty$ and all the volumes in this article are counted with multiplicity). When $V= \Cc$ and $\tau= \| z
\|^2$, M. McQuillan constructed such currents in \cite{McQ} (see
\cite{Br} too). These currents are fundamental tools in the study of the
hyperbolicity of $X$ (see for example \cite{Du1}). When
the dimension of $V$ is higher or equal to $2$ it is not always
possible to produce Ahlfors' currents. Indeed, for example, there exists domains
$\Omega$ in $\Cc^2$ which are biholomorphic to $\Cc^2$ and such that
$\overline{\Omega} \neq \Cc^2$ (Fatou-Bieberbach domains). As a consequence, to produce
Ahlfors' currents it is necessary to add a hypothesis on $f$.

When the dimension of $X$ is equal to $k$, there exists some criteria
which imply that $f(V)$ is dense in $X$ (see \cite{Ch}, \cite{St},
\cite{Wu}, \cite{Hi}, \cite{Gr}, \cite{CG} and \cite{Si}). These criteria used the
degrees of $f$ (see \cite{Ch}) or the growth of
the function $f$.

Our goal is to give criteria which use these degrees in order to produce Ahlfors'
currents in $X$. Of course, in the case where the dimension of $X$ is
equal to $k$, the existence of such currents will automatically
imply that $f(V)$ is dense in $X$. Indeed, $[X]$ is the only
positive closed current of bidimension $(k,k)$ in $X$ (modulo a
normalization).

In this article, we will use the following degrees ($t_{k-1}$ will be slightly different
from the Chern's one):

$$t_k(r)= \int_{V(r)} f^{*} \omega^k,$$

$$t_{k-1}(r)= \int_{V(r)}  i \partial \tau \wedge
\overline{\partial \tau}  \wedge f^{*} \omega^{k-1}           .$$

Let $\mathcal{C}$ be the set of critical values of $\tau$ in $[r_0 ,
  +\infty[$. V is connected and non-compact so we can suppose that $[r_0 , + \infty
  [  \subset \tau(V)$.

The criteria that we will give on $t_k$ and $t_{k-1}$ will strongly use the
following  inequality:

\begin{theorem}{\label{inequality}}

The functions $t_k$ and $t_{k-1}$ are $C^1$ on $]r_0 , +
  \infty[ \setminus \mathcal{C}$ and $C^0$ on $]r_0 , +
  \infty[$. If $r \in  ]r_0 , +
  \infty[ \setminus \mathcal{C}$ then

$$\| \partial f_{*} [V(r) ] \|^2 \leq K(X) t'_{k-1}(r) t'_k(r).$$

Here $\|.\|$ is the norm in the sense of currents and $K(X)$ is a
constant which depends only on $X$.

\end{theorem}

By using the previous inequality we can prove some criteria which imply the
existence of Ahlfors' currents. Here we give the two following ones:

\begin{theorem}{\label{criterion1}}

We suppose that $f$ is nondegenerate and of finite type (i.e. there exists $C_1
\mbox{, } C_2 \mbox{, } r_1 >0$ such that $ \mbox{volume}(f(V(r)))
\leq C_1 r^{C_2}$ for $r \geq r_1$).

If 

$$ \limsup_{r \rightarrow + \infty} \frac{t_{k-1}(r)}{r^2 t_k(r)}=0$$

then there exists a sequence $r_n$ which goes to infinity such that
$\frac{f_{*} [V(r_n)]}{\mbox{volume}(f(V(r_n)))}$ converges to a closed
positive current with bidimension $(k,k)$ and mass equal to $1$.

\end{theorem}

When $V=\Cc$ and $\tau= \| z \|^2$, the hypothesis of finite type is
true modulo a Brody's renormalization (see for example \cite{La}).

We give now one criterion which doesn't use this hypothesis.

\begin{theorem}{\label{criterion2}}

If $f$ is nondegenerate and if there exists $\epsilon >0$ and $L>0$
such that:

$$ \limsup_{ r \notin \mathcal{C} \mbox{, } r \rightarrow + \infty} \frac{t'_{k-1}(r)}{r t_k^{1-
    \epsilon}(r)} \leq L$$

then there exists a sequence $r_n$ which goes to infinity such that
$\frac{f_{*} [V(r_n)]}{\mbox{volume}(f(V(r_n)))}$ converges to a closed
positive current with bidimension $(k,k)$ and mass equal to $1$.

\end{theorem}

The plan of this article is the following one: in the first part we prove
the inequality (theorem \ref{inequality}), in the second one we give
the proof of the both criteria (theorems \ref{criterion1} and \ref{criterion2}). In the third part, we will
give a new formulation of the criteria for the
particular case where $V=\Cc^k$.

\section{\bf Proof of the inequality}

Let $\mathcal{C}$ be the set of critical values of $\tau$ in $[r_0 ,
  +\infty[$. We recall that we can suppose $[r_0 , + \infty
  [  \subset \tau(V)$. Notice that the point (iii) of the hypothesis on $\tau$ implies that $\mathcal{C}$ is
  discrete. When $r \in ]r_0 , + \infty[$ and $r \notin \mathcal{C}$
  then $\tau : \tau^{-1}( ]r- \epsilon, r+ \epsilon[) \mapsto ]r-
  \epsilon, r+ \epsilon[$ is a submersion for $\epsilon >0$ small enough. In particular,
  $\tau^{-1}(r)$ is a submanifold of $V$ and $\partial V(r)=
  \tau^{-1}(r)$. When $r \in \mathcal{C}$, then $\tau^{-1}(r)$ is a
  compact set which is a submanifold of $V$ outside a neighbourhood
  of a finite number of points.

We begin now with the following lemma:

\begin{lemma}{\label{lemma1}}

The functions $t_k$ and $t_{k-1}$ are $C^1$ on $]r_0 , +
  \infty[ \setminus \mathcal{C}$ and $C^0$ on $]r_0 , +
  \infty[$.

\end{lemma}

\begin{proof}

The form $f^{*} \omega^k$ is positive and smooth and $i \partial \tau
\wedge \overline{\partial \tau} \wedge f^{*} \omega^{k-1}$ is positive and continuous ($\tau$ is $C^1$) so it
is enough to show that $t(r)=\int_{V(r)} \Phi$ is $C^1$ on
$]r_0 , +  \infty[ \setminus \mathcal{C}$ and $C^0$ on $]r_0 , +
  \infty[$ with $\Phi$ a positive continuous
    form of bidegree $(k,k)$.

We take $r \in ]r_0 , +  \infty[ \setminus \mathcal{C}$ and $\epsilon >0$
  such that $\tau : \tau^{-1}( ]r- \epsilon, r+ \epsilon[) \mapsto ]r-
  \epsilon, r+ \epsilon[$ is a submersion.  Now, if $r' \in
    ]r- \epsilon, r[$, we have:

$$\frac{t(r)-t(r')}{r-r'}= \frac{1}{r-r'} \int_{\tau^{-1}([r',r[)}
	  \Phi=  \frac{1}{r-r'} \int_{[r',r[} \tau_{*} \Phi.$$

The form $\tau_{*} \Phi$ is continuous so it is equal to $\alpha(s)
ds$ with $\alpha$ in $C^0 (]r- \epsilon, r+ \epsilon[)$. We obtain:

$$\frac{t(r)-t(r')}{r-r'}=\frac{1}{r-r'} \int_{r'}^{r} \alpha(s) ds$$

which converges to $\alpha(r)$ when $r' \rightarrow r$. The same thing
happens when we consider $r' \in
    ]r, r+ \epsilon[$, so the function $t$ is differentiable at $r$
	      and $t'(r)=\alpha(r)$. In particular $t$ is $C^1$ on $]r_0 , +
  \infty[ \setminus \mathcal{C}$.

\begin{remark}{\label{rem1}}

Notice that here we did not use that $\Phi$ is positive. We will use this remark in the proof of the theorem \ref{inequality}.

\end{remark}

Now, consider $r \in \mathcal{C}$. If we take $\epsilon > 0$, then we
can find two neighbourhoods $W_{\epsilon} \Subset W_{2 \epsilon} $ of the (finite) number of the
critical points in $\{ \tau=r \}$ such that $\int_{W_{2 \epsilon}} \Phi \leq
\epsilon$ (because $\Phi$ is continuous). Now, let $\psi$ be a
$C^{\infty}$ function which is equal to $1$ in a neighbourhood of $\overline{W_{\epsilon}}$ and to
$0$ outside $W_{2 \epsilon}$ ($0 \leq \psi \leq 1$). Then, if $r'<r$,

$$t(r)-t(r')= \int_{V(r) \setminus V(r')} \psi \Phi +  \int_{V(r)
  \setminus V(r')} (1 - \psi) \Phi \leq \epsilon + \int_{V(r)
  \setminus V(r')} ( 1 - \psi) \Phi.$$

If $\alpha>0$ is small then $\tau$ is a submersion on $\tau^{-1}(]r-
  \alpha, r + \alpha[) \cap (V \setminus W_{\epsilon})$. In
  particular the function 
$$r' \mapsto  \int_{V(r)  \setminus V(r')} (1- \psi) \Phi=
  \int_{r'}^{r} \tau_{*}((1-\psi) \Phi)$$
goes to $0$ when $r' \rightarrow
  r$. The same thing happens
  when we take $r' > r$. As a consequence, there exists $\delta >0$ such that if $|r - r'|
  < \delta$ then $| t(r)-t(r') | \leq 2 \epsilon$, i.e. $t$ is
  continuous at $r$.

\end{proof}

We give now the proof of the theorem \ref{inequality}.

We take $r \in ]r_0 , +  \infty[ \setminus \mathcal{C}$. We have:

$$\| \partial f_{*} [V(r)] \| = \sup_{\Psi \in \mathcal{F}(k-1,k)}
    \left| \langle \partial f_{*} [V(r)] , \Psi \rangle \right|$$

where $\mathcal{F}(k-1,k)$ is the set of smooth $(k-1,k)$ forms $\Psi$
with $\| \Psi \| \leq 1$. If $\Psi \in  \mathcal{F}(k-1,k)$ then we
can write

$$\Psi= \sum_{i=1}^{K(X)} \theta_i \wedge \Omega_i$$

where $K(X)$ is a constant which depends only on $X$, the $\theta_i$
are smooth forms of bidegree $(0,1)$ with $\| \theta_i \| \leq1$ and the $\Omega_i$
are (strongly) positive smooth forms of bidegree $(k-1,k-1)$ with $\| \Omega_i \|
\leq K(X)$. So, to prove the inequality it is
sufficient to bound from above $ | \langle \partial f_{*} [V(r)] ,
\theta \wedge \Omega \rangle |^2$ by $ K(X) t'_{k-1}(r) t'_k(r)$ with $\theta$ a smooth form of bidegree $(0,1)$ with
$\| \theta \| \leq 1$ and $\Omega$ a positive smooth form of
bidegree $(k-1,k-1)$ with $\| \Omega \| \leq 1$.

If $\epsilon >0$ is small then $\tau : \tau^{-1}( ]r- \epsilon, r+ \epsilon[) \mapsto ]r-
  \epsilon, r+ \epsilon[$ is a submersion. Now, if we take $r' \in ]r- \epsilon, r[$, we have:

$$A(r',r):=\left| \frac{1}{r-r'} \int_{r'}^{r} \langle  \partial
	      f_{*} [V(s)] ,   \theta \wedge \Omega \rangle ds \right|= \left|  \frac{1}{r-r'}
	  \int_{r'}^{r} \langle  \partial  [V(s)]  , f^{*} \theta
	      \wedge f^{*} \Omega  \rangle ds \right|.  $$

If we use the Stokes' theorem, we have:
 
$$A(r',r)= \left|  \frac{1}{r-r'}
	  \int_{r'}^{r} \langle   [\partial V(s)]  , f^{*} \theta
	      \wedge f^{*} \Omega  \rangle ds \right|=  \left|
	  \frac{1}{r-r'} \int_{r'}^{r} \left\langle  [\tau=s] , f^{*}
\theta \wedge f^{*} \Omega \right\rangle ds \right|,$$

because for $s \in ]r- \epsilon, r+ \epsilon[$ the boundary of $V(s)$
    is $\{ \tau = s \}$.

We obtain:

$$A(r',r)= \left| \frac{1}{r-r'} \int_{r'}^{r} \left( \int_{\tau=s}  f^{*}
\theta \wedge f^{*} \Omega \right) ds \right|.$$

 Now $\tau : \tau^{-1}( ]r- \epsilon, r+ \epsilon[) \mapsto ]r-
  \epsilon, r+ \epsilon[$ is a submersion, so by using Fubini's
    theorem  (see \cite{Chi} p. 334), we have:

$$A(r',r)= \left| \frac{1}{r-r'} \int_{V(r) \setminus V(r')} d \tau
    \wedge  f^{*} \theta \wedge f^{*} \Omega \right| =   \left| \frac{1}{r-r'} \int_{V(r) \setminus V(r')} \partial \tau
    \wedge  f^{*} \theta \wedge f^{*} \Omega \right|.$$

Now, if we consider,

$$\{ \phi, \psi \} := \int_{V(r) \setminus V(r')} i \phi \wedge \overline{\psi} \wedge f^{*} \Omega$$

where $\phi$ and $\psi$ are continuous forms of bidegree $(1,0)$, then
$\{ \phi, \phi \} \geq 0$ (because $\Omega$ is positive) and so by using the proof of the
Cauchy-Schwarz's inequality we obtain that:

$$| \{ \phi, \psi \} | \leq (\{ \phi, \phi \} )^{1/2} (\{ \psi, \psi
\} )^{1/2}.$$

In particular,

$$A(r',r)^2 \leq \left|  \frac{1}{r-r'}
	      \int_{V(r) \setminus V(r')} i \partial \tau \wedge
	      \overline{\partial \tau} \wedge f^{*} \Omega \right|
	      \times \left|  \frac{1}{r-r'}
	      \int_{V(r) \setminus V(r')} i  \overline{f^{*} \theta}
	      \wedge  f^{*} \theta \wedge f^{*} \Omega \right|.$$

Now $i  \overline{f^{*} \theta}  \wedge  f^{*} \theta \wedge f^{*}
\Omega$ is equal to $f^{*}(i \overline{\theta} \wedge \theta \wedge
\Omega)$ and $i \overline{\theta} \wedge \theta \wedge \Omega \leq K'(X)
\omega^k$ (which means that $K'(X)\omega^k - i \overline{\theta} \wedge
\theta \wedge \Omega$ is a (strongly) positive form). Here $K'(X)$
depends only on $X$ because $\| \theta \| \leq 1$ and $\| \Omega \|
\leq 1$.

As a consequence, we have:

$$ \left|  \frac{1}{r-r'}
	      \int_{V(r) \setminus V(r')} i  \overline{f^{*} \theta}
	      \wedge  f^{*} \theta \wedge f^{*} \Omega \right| \leq
	      K'(X) \left|  \frac{1}{r-r'}
	      \int_{V(r) \setminus V(r')} f^{*} \omega^k \right|=
	      K'(X) \left( \frac{t_k(r)-t_k(r')}{r-r'} \right).$$

On the other hand, there exists a constant $K''(X)$ with $\Omega \leq K''(X)
\omega^{k-1}$ (we use $\| \Omega \| \leq 1$). So, we have

$$\left|  \frac{1}{r-r'}
	      \int_{V(r) \setminus V(r')} i \partial \tau \wedge
	      \overline{\partial \tau} \wedge f^{*} \Omega \right|
	      \leq K''(X) \left( \frac{t_{k-1}(r)-t_{k-1}(r')}{r-r'} \right).$$

We obtain:

\begin{equation}{\label{eq1}}
A(r',r)^2 \leq K(X) \left( \frac{t_{k-1}(r)-t_{k-1}(r')}{r-r'}
\right) \left( \frac{t_k(r)-t_k(r')}{r-r'} \right).\\
\end{equation}

Now, when $r' \rightarrow r$ 
$$A(r',r)^2 \rightarrow \left| \langle \partial f_{*} [V(r)] , \theta
\wedge \Omega \rangle \right|^2$$

because the function $s \mapsto \langle  \partial f_{*} [V(s)] ,
\theta \wedge \Omega \rangle= \int_{V(s)} \partial f^{*} (\theta \wedge \Omega)$ is continuous on $]r- \epsilon, r+ \epsilon[$ (see remark \ref{rem1}).

Finally, if we take $r' \rightarrow r$ in the inequality (\ref{eq1}), we have:

$$ \left| \langle \partial f_{*} [V(r)] , \theta \wedge \Omega \rangle \right|^2 \leq
K(X) t'_{k-1}(r) t'_k(r) $$

which gives the wanted inequality.
 
\section{\bf Proof of the theorems \ref{criterion1} and \ref{criterion2}}

\subsection{{\bf Proof of the first criterion}}

We begin with this lemma:

\begin{lemma}

If $f$ is nondegenerate and of finite type then there exits a
constant $K>0$ such that:

$$\forall r_2 > 0 \mbox{  } \exists r \geq r_2 \mbox{ with }
\mbox{volume}(f(V(2r))) \leq K \mbox{volume}(f(V(r))).$$

\end{lemma}

\begin{proof}

The hypothesis implies that there exists $C_1
\mbox{, } C_2 \mbox{, } r_1 >0$ such that $ \mbox{volume}(f(V(r)))
\leq C_1 r^{C_2}$ for $r \geq r_1$.

If the conclusion of the lemma is false then for all $K>0$ there
exists $r_2 >0$ such that for all $r \geq r_2$ we have
$\mbox{volume}(f(V(2r))) \geq K \mbox{volume}(f(V(r)))$.

So, if we take $K >> 2^{C_2}$ then we obtain (if $l$ is high enough):

$$C_1(2^lr_2)^{C_2} \geq \mbox{volume}(f(V(2^lr_2))) \geq K^l \mbox{volume}(f(V(r_2))).$$

As a consequence we have

$$\mbox{volume}(f(V(r_2))) \leq C_1r_2^{C_2} \left( \frac{2^{C_2}}{K}
\right)^l$$

which implies that $\mbox{volume}(f(V(r_2)))=0$ when we take $l
\rightarrow \infty$. It contradicts the fact that $f$ is nondegenerate.

\end{proof}

By using this lemma, we can find a sequence $R_n \rightarrow + \infty$
which satisfies

$$\mbox{volume}(f(V(2R_n))) \leq K \mbox{volume}(f(V(R_n))).$$

The theorem \ref{inequality} gives now that:

$$\int_{R_n}^{2R_n} \| \partial f_{*} [V(r)] \| dr \leq K(X)
\int_{R_n}^{2R_n} \sqrt{t'_{k-1}(r)} \sqrt{t'_k(r)} dr.$$

We use the following sense to the integrals: for example, if
there is one point $a_n$ of $\mathcal{C}$ in $[R_n , 2 R_n]$, we
consider $\int_{R_n}^{2R_n} = \lim_{\epsilon \rightarrow 0} 
\int_{[R_n, a_n - \epsilon] \cup [a_n+  \epsilon, 2R_n]} $. All the functions that we consider are non
negative, so the limit exists in $[0 , + \infty]$. 

Now, by using the Cauchy-Schwarz's inequality, the last integral is smaller than

$$K(X) \left( \int_{R_n}^{2R_n} t'_{k-1}(r) dr \right)^{1/2}  \left(
\int_{R_n}^{2R_n} t'_{k}(r) dr \right)^{1/2} \leq K(X)
\sqrt{t_{k-1}(2R_n)}  \sqrt{t_{k}(2R_n)}.$$

For the last inequality it is important to use that $t_{k-1}$ and $t_k$
are continuous on $]r_0, + \infty[$ (see theorem \ref{inequality}).

It implies that there exists a sequence $r_n \in [R_n , 2R_n]$ such
that:

$$\| \partial f_{*} [V(r_n)] \| \leq \frac{K(X)}{R_n}
\sqrt{t_{k-1}(2R_n)}  \sqrt{t_{k}(2R_n)},$$

i.e.

$$\frac{\| \partial f_{*} [V(r_n)] \|}{\mbox{volume}(f(V(r_n)))} \leq 2K(X)
\sqrt{\frac{t_{k-1}(2R_n)}{(2R_n)^2 t_k(2R_n)}}
\frac{t_k(2R_n)}{t_k(r_n)}$$

because $\mbox{volume}(f(V(r_n)))=t_k(r_n)$.

Now, we have 

$$\frac{t_k(2R_n)}{t_k(r_n)} \leq \frac{t_k(2R_n)}{t_k(R_n)} \leq K$$

and by using the hypothesis,

$$\sqrt{\frac{t_{k-1}(2R_n)}{(2R_n)^2 t_k(2R_n)}} \rightarrow 0.$$

So, we obtain that 

$$\frac{\| \partial f_{*} [V(r_n)] \|}{\mbox{volume}(f(V(r_n)))}
\rightarrow 0.$$

The current $T_n:=\frac{f_{*} [V(r_n)]}{\mbox{volume}(f(V(r_n)))}$ is positive with bidimension $(k,k)$ and mass equal to $1$, so there exits a
subsequence of $(T_n)$ which converges to a positive current $T$ with bidimension
$(k,k)$ and mass $1$. Moreover, 

$$\| \partial T_n \| = \frac{\| \partial f_{*} [V(r_n)] \|}{\mbox{volume}(f(V(r_n)))}
\rightarrow 0,$$

so the limit current $T$ is closed. It proves the first criterion.

\subsection{{\bf Proof of the second criterion}}

Take $\epsilon>0$ and $L>0$ such that

$$ \limsup_{r \notin \mathcal{C} \mbox{, }r \rightarrow + \infty} \frac{t'_{k-1}(r)}{r t_k(r)^{1-
    \epsilon}} \leq L.$$

Let $R_n$ be a sequence of positive reals which goes to $+
\infty$. By using the theorem \ref{inequality}, we have (see the proof
of the last criterion for the definition of the integrals):

$$\int_{r_0 +1}^{R_n} \frac{\| \partial f_{*} [V(r)] \|^2}{t'_{k-1}(r)
  t_k(r)^{1+ \epsilon}} dr \leq K(X) \int_{r_0 +1}^{ R_n}
  \frac{t'_k(r)}{t_k(r)^{1+ \epsilon}} dr.$$

This last integral is smaller than $\frac{K(X)}{\epsilon
  t_k(r_0 +1)^{\epsilon}} \leq K'(X,f)$ (here we use the fact that $\frac{1}{t_k(r)}$
  is continuous on $]r_0 , + \infty[$).

So, we have 

$$\int_{r_0 +1}^{+ \infty} \frac{1}{r} \left( \frac{r\| \partial f_{*}
  [V(r)] \|^2}{t'_{k-1}(r)
  t_k(r)^{1+ \epsilon}} \right) dr \leq K'(X,f),$$

and $\int_{r_0 +1}^{+ \infty} \frac{1}{r} dr = + \infty$ implies that there exits a sequence $r_n \rightarrow + \infty$ such
that $r_n \notin \mathcal{C}$ and:

$$\epsilon(n):=\frac{r_n \| \partial f_{*} [V(r_n)] \|^2}{t'_{k-1}(r_n)
  t_k(r_n)^{1+ \epsilon}} \rightarrow 0.$$

We obtain

$$ \left( \frac{\| \partial f_{*} [V(r_n)] \|}{\mbox{volume}(f(V(r_n)))}
\right)^2 = \frac{\epsilon(n)}{r_n}  \frac{ t'_{k-1}(r_n)}{ t_k(r_n)^{1 - \epsilon}} \leq (L+1) \epsilon(n),$$

by hypothesis (for $n$ high enough).

So, 

$$\frac{\| \partial f_{*} [V(r_n)] \|}{\mbox{volume}(f(V(r_n)))}
\rightarrow 0.$$

Now, by using exactly the same argument than in the proof of the
previous criterion, we obtain that there exists a subsequence of
$T_n:=\frac{f_{*} [V(r_n)]}{\mbox{volume}(f(V(r_n)))}$ which converges to
a closed positive current of bidimension $(k,k)$ and with mass equal
to $1$.

\section{\bf{The particular case $V=\Cc^k$}}

In this paragraph we consider the particular case where $V= \Cc^k$. 

Let $\beta$ be the standard kähler form in $\Cc^k$. We want to transform
our previous criteria by using $\beta$ instead of $i \partial \tau \wedge \overline{\partial \tau}$. More precisely, we consider:

$$a_k(r)= \int_{B(0,r)} f^{*} \omega^k$$

and

$$a_{k-1}(r)= \int_{B(0,r)} \beta \wedge  f^{*} \omega^{k-1}.$$

Then we can prove a new formulation of our three theorems:

\begin{theorem}{\label{inequality2}}

The functions $a_k$ and $a_{k-1}$ are $C^1$ on $]0 , +
  \infty[$ and for $r >0$ we have

$$\| \partial f_{*} [B(0,r) ] \|^2 \leq K(X) a'_{k-1}(r) a'_k(r).$$

Here $\|.\|$ is the norm in the sense of currents and $K(X)$ is a
constant which depends only on $X$.
 
\end{theorem}

\begin{proof}

We apply the theorem \ref{inequality} with $V= \Cc^k$ and $\tau= \|
z \|^2$ (here we have $\mathcal{C} = \{ 0 \}$) and then for $r>0$:

$$\| \partial f_{*} [V(r^2) ] \|^2 \leq K'(X) t'_{k-1}(r^2) t'_k(r^2).$$

Now, $a_k(r)=t_k(r^2)$, so $a_k$ is $C^1$ in $]0 , +  \infty[$ and 

$$t'_k(r^2) = \frac{a'_k(r)}{2r}.$$

The function $a_{k-1}(r)=t(r^2)$ with $t(r)=\int_{V(r)} \beta \wedge f^{*}
\omega^{k-1} $ so $a_{k-1}$ is $C^1$ in $]0 , +  \infty[$
    (see the proof of the lemma \ref{lemma1}).

Moreover, 

$$t_{k-1}(r^2)= \int_{V(r^2)}  i \partial \tau \wedge
\overline{\partial \tau} \wedge  f^{*} \omega^{k-1}             =
\int_{B(0,r)}  i \partial \tau \wedge \overline{\partial \tau} \wedge
f^{*} \omega^{k-1} ,$$

and $i \partial \tau \wedge \overline{\partial \tau}=i \sum_{i,j} \overline{z_i}
z_j d z_i \wedge d \overline{z_j}$.

On $B(0,r)$ this last form is smaller than $K(k) \beta r^2$.

If we take $0<r'<r$ then

$$t_{k-1}(r^2) - t_{k-1}(r'^2)= \int_{B(0,r) \setminus B(0,r')}  i
    \partial \tau \wedge \overline{\partial \tau} \wedge
    f^{*} \omega^{k-1}  \leq K(k) r^2 \int_{B(0,r) \setminus B(0,r')}
    \beta \wedge f^{*} \omega^{k-1} .$$

If we divide by $r-r'$ and if we take the limit $r' \rightarrow r$
then we obtain:

$$2r t'_{k-1}(r^2) \leq K(k) r^2 a'_{k-1}(r).$$

Finally, we have:

$$\| \partial f_{*} [B(0,r) ] \|^2= \| \partial f_{*} [V(r^2) ] \|^2  \leq
K'(X) t'_{k-1}(r^2) t'_k(r^2) \leq K(X) a'_{k-1}(r) a'_k(r),$$

with $K(X)=K(k)K'(X)$ (we recall that the dimension of $X$ is higher
or equal to $k$). This is the inequality that we wanted.

\end{proof}

Now if we remplace in the proof of the theorems \ref{criterion1} and
\ref{criterion2} the function $t_{k-1}$ by $a_{k-1}$, the function $t_k$ by $a_k$ and $V(r)$ by
$B(0,r)$ then we obtain the two following criteria:

\begin{theorem}

We suppose that $f$ is nondegenerate and with finite type (i.e. there exists $C_1
\mbox{, } C_2 \mbox{, } r_1 >0$ such that $ \mbox{volume}(f(B(0,r)))
\leq C_1 r^{C_2}$ for $r \geq r_1$).

If 

$$ \limsup_{r \rightarrow + \infty} \frac{a_{k-1}(r)}{r^2 a_k(r)}=0$$

then there exists a sequence $r_n$ which goes to infinity such that
$\frac{f_{*} [B(0,r_n)]}{\mbox{volume}(f(B(0,r_n)))}$ converges to a closed
positive current with bidimension $(k,k)$ and mass equal to $1$.

\end{theorem}

\begin{theorem}

If $f$ is nondegenerate and if there exists $\epsilon >0$ and $L >0$
such that:

$$ \limsup_{r \rightarrow + \infty} \frac{a'_{k-1}(r)}{r a_k^{1-
    \epsilon}(r)} \leq L$$

then there exists a sequence $r_n$ which goes to infinity such that
$\frac{f_{*} [B(0,r_n)]}{\mbox{volume}(f(B(0,r_n)))}$ converges to a closed
positive current with bidimension $(k,k)$ and mass equal to $1$.

\end{theorem}

Notice that when $k=1$ then $a_{k-1}(r)=\pi r^2$ and so, in this
context, the hypothesis
of this criterion is always true if $f$ is nondegenerate.

\bigskip

Henry de Thélin\\
Université Paris-Sud (Paris 11)\\
Mathématique, Bât. 425\\
91405 Orsay\\
France


\begin{thebibliography}{00}

\bibitem{Br} M. Brunella, \textit{Courbes entières et feuilletages
  holomorphes}, Enseign. Math., {\bf 45} (1999), 195-216.

\bibitem{CG} J.A. Carlson and P. Griffiths, \textit{The order
  functions for entire holomorphic mappings}, Proc. Tulane Univ. Program, (1974), 225-248.

\bibitem{Ch} S.-S. Chern, \textit{The integrated form of the first main
theorem for complex analytic mappings in several complex variables},
  Ann. of Math. (2), {\bf 71} (1960), 536-551.

\bibitem{Chi} E.M. Chirka, \textit{Complex analytic sets}, Kluwer
  Academic Publishers (1989).

\bibitem{Du1} J. Duval, \textit{Singularités des courants d'Ahlfors},
  Ann. Sci. Ecole Norm. Sup., {\bf 39} (2006), 527-533.

\bibitem{Gr} P. Griffiths, \textit{Some remarks on Nevanlinna
  theory}, Proc. Tulane Univ. Program, (1974), 1-11.

\bibitem{Hi} J.J. Hirschfelder, \textit{The first main theorem of
  value distribution in several variables}, Invent. Math., {\bf 8}
  (1969), 1-33.

\bibitem{La} S. Lang, \textit{Introduction to complex hyperbolic
  spaces}, Springer-Verlag (1987).

\bibitem{McQ} M. McQuillan, \textit{Diophantine approximations and
  foliations}, Inst. Hautes Etudes Sci. Publ. Math., {\bf 87} (1998), 121-174.

\bibitem{Ra} R.M. Range, \textit{Holomorphic functions and integral
  representations in several complex variables}, Springer-Verlag
  (1986).

\bibitem{Si} N. Sibony and P.M. Wong, \textit{Some remarks on the
  Casorati-Weierstrass theorem}, Ann. Polon. Math., {\bf 39} (1981), 165-174.


\bibitem{St} W. Stoll, \textit{A general first main theorem of value
  distribution}, Acta Math., {\bf 118} (1967), 111-191.

\bibitem{Wu} H. Wu, \textit{Remarks on the first main theorem in
  equidistribution theory II}, J. Differential Geometry, {\bf 2} (1968), 369-384.




\end{thebibliography}
\end{document}